\documentstyle{article}

\def \l {\lambda}

\def \E {{\bf E}}
\def \s {\sigma}

\def \N {{\bf N}}

\def \V {\Vert}
\def \v {\vert}

\def \la {\langle}
\def \ra {\rangle}

\def \T {{\hbox{ Tr }}}

\def \d {{\hbox {d}}}
\def \vep {\varepsilon}
\def \a {\alpha}

\def \g {\gamma}

\def \l {\lambda}

\def \CN {{\cal N}}

\begin{large}
\author{A.Khorunzhy\thanks{The author 
thanks the Department of Mathematics and Statistics 
of the Brunel University of West London, where 
this work has been started. The financial support of the 
Royal Society (London) is gratefully acknowledged. }\\
Laboratoire de Physique Math\'ematique et G\'eometrie\\
Universit\'e Paris 7 - Denis Diderot, FRANCE\thanks{On leave of absence
from Mathematical Division, Institute for Low Temperature Physics, Kharkov, UKRAINE}}

\title{SPARSE RANDOM MATRICES AND STATISTICS OF ROOTED TREES}

\begin{document}

\maketitle
\begin{abstract}
We consider the ensemble of $N\times N$ random symmetric matrices
$\{A_{N,p}\}$ with independent entries such that 
$A_{N,p}$ has, in average, $p$ non-zero elements per row.
We study the asymptotic behaviour of the maximal (in magnitude)
eigenvalue
$\l_{\max}(N)$ of matrices $A_{N,p}/\sqrt p$ when they are 
large and sparse, i.e. in the limit $N,p\to\infty,$
$p=o(N)$. We prove that the value 
$p_c= \log N$ is the critical one for $\l_{\max}(\infty)$
to be bounded.

Our arguments are based on the  calculus of the tree-type
graphs for  random matrices with independent entries. Asymptotic
properties of sparse  random matrices essentially
depend on the typical degree of a tree vertex that we prove to be
finite.
\end{abstract}
\vskip 0.2cm
{\it AMS Classification numbers (1991):} 15A52 (Primary) 05C30, 60F99 (secondary) 
\vskip 0.2cm
{\it Key words:} random matrices; eigenvalue distribution;
enumeration of graphs; Erd\"os-R\'enyi partial sums.

\section{Introduction}The study of eigenvalue distribution of 
$N$-dimensional random matrices in the limit 
$N\to\infty$ has been initiated by E.Wigner \cite{W}.
He considered the ensemble of real symmetric matrices 
$$
A_N(x,y) = a(x,y),\quad x,y = 1,\dots, N
\eqno (1.1)
$$
where $\{a(x,y),x\le y, \, x,y \in \N\}$ are jointly independent random variables
with zero mean and variance $1$. The celebrated semicircle
(or Wigner) law states that the  distribution function 
of eigenvalues  $\l_1^{(N)}\le\dots\le\l_N^{(N)}$ of the matrix
$\hat A_N = A_N/\sqrt N$ given by relation
$$
\s(\l;\hat A_N) = \#\{j: \, \l_j^{(N)}\le \l\} N^{-1}
\eqno (1.2)
$$
converges in the average as $N\to\infty$ to 
the nonrandom distribution $\s_w(\l)$
$$
\s'_w(\l) = \cases{(2\pi)^{-1}\sqrt{4-\l^2},& if $\vert \l\vert\le 2$ \cr
                             0, & if $\vert \l\vert>2$\cr}
\eqno (1.3)
$$
provided all even moments of $a(x,y)$ exist and the odd ones vanish.

Wigner introduced random matrices (1.1) to model distribution of energy levels
of heavy atomic nuclei. Such a nucleus consists of large number
($N\sim 100$) of particles that interact each with other. According to this,
the statistical description naturally presumes consideration 
of random matrices whose entries are random variables of the same order
of magnitude (see e.g. \cite{P}).

Recent applications of random matrices lead to various generalizations
of the Wigner ensemble $\{\hat A_N\}$. In particular, in the neural network
theory the dilute version of (1.1) is used to describe the system, where some
randomly chosen links between elements are broken 
(see e.g. \cite{A}). Such a dilute matrix can be determined as 
$$
A_{N,d}(x,y) = a(x,y) d_N(x,y), \quad d_N(x,y)= d_N(y,x),
$$
where $\Lambda_N\equiv \{ d_N(x,y),\,x\le y,\,x,y =1,\dots,N\}$ 
is a family of independent random
variables, also independent from $\{a(x,y)\}$, with the distribution
$$
d_N(x,y) = 
 \cases{1,& with probability $p/N$ ;\cr
                             0, & with probability $1-p/N$;\cr}
                             \quad p\le N.
\eqno (1.4)
$$
It is not hard to see that the  $A_{N,d}$ has, in average,
$p$ non-zero elements per row. Thus,  $\{A_{N,d}\}$ could be regarded as the
ensemble of matrices intermediate between the discrete version of
random Schroedinger operator and the Wigner random matrices.
This explains the interest to dilute random matrices from the spectral theory
point of view (see e.g. \cite{KKPS}).

The asymptotic properties of the eigenvalue distribution of 
$A_N^{(d)}$
were considered in different aspects \cite{E,MF,RB}. In particular, 
it was shown in \cite{RB} and proved in 
\cite{KKPS} that the eigenvalue distribution
of the matrix
$$
 A_N^{(p)}(x,y) = a(x,y) \hat d_N(x,y), \quad 
\hat d_N(x,y)= {1\over \sqrt p} d_N(x,y)
\eqno (1.5)
$$
converges in probability in the limit $N,p\to \infty$ to the semicircle
distribution $\s_w(\l)$ (1.3);
$$
{\hbox{p-}}\lim_{N,p\to\infty} \s(\l;A_N^{(p)})= \s_w(\l)
\eqno (1.6)
$$
that means that the semicircle law holds for dilute random matrices.

The eigenvalue distribution function concerns a fraction of eigenvalues
of the matrix. In applications 
it is often important to know whether the norm of 
$A_{N}^{(p)}$ is bounded (see e.g. \cite{CV}), i.e. whether there 
are eigenvalues of $A_{N}^{(p)}$ outside of the support of $\s_w'(\l)$
in the limit $N,p\to\infty$.

In present paper we consider the
probability distribution of the
spectral norm
$$
\V  A_N^{(p)}\V = \max\{\v \l_1^{(N,p)}\v,\,\v\l_N^{(N,p)}\v\}\equiv
\l_{\max}^{(N,p)},
$$
where 
$\l_1^{(N,p)}\le\dots\le\l_N^{(N,p)}$  are eigenvalues of 
$ A_N^{(p)}$.
We study the asymptotic behaviour of 
$\V  A_N^{(p)}\V$ 
when matrices
$ A_N^{(p)}$ are large and sparse, i.e. in the limit
$N,p\to\infty,\, p=o(N)$. Our results show that the rate 
$p_c=\log N$ is the critical one for the limit of $\V  A_N^{(p)}\V$
to be either bounded (in this case it is equal to 2)
or unbounded.

\section{Main results and scheme of the proof}

First of all, let us note
that one can determine 
the collection of random variables $\Lambda\equiv \{\Lambda_N,\, N\in \N\}$ on the same
probability space
$\Omega_d$
with the help of the procedure due to C.Newman (see e.g. \cite{BG}).
We also determine $\{a(x,y),\, x,y \in \N\}$ on the same $\Omega_a$.

We denote the measure and mathematical expectation corresponding to $\Lambda$ by $\mu_d$
and  angles $\la\cdot\ra$, respectively.
We also denote by $\E\{\cdot\}$ the mathematical expectation with respect to the
measure $\mu_a$ generated by the random variables $\{a(x,y),\, x,y \in \N\}$.

We assume that the random variables 
$a(x,y)$ satisfy conditions:
$$
\E [a(x,y)]^{2k+1} = 0  \quad \forall x,y,k\in\N,
\eqno (2.1a)
$$
and
$$
\E [a(x,y)]^{2}= 1, \quad \sup_{x,y}\E [a(x,y)]^{2k} = V_{2k} <\infty
\quad
\forall k\in \N.
\eqno (2.1b)
$$

{\bf Theorem 2.1}

{\it Assume that the random variables $a(x,y)$ 
are such that there exists $\delta >0$ such that 
$$
V_{2k}\le k^{\delta k} \quad \forall k\in \N.
\eqno (2.2)
$$
Then in the limit $N,p\to \infty$ such that  $p$ satisfies condition
$$
 {p\over (\log N)^{1+\delta'}}\to \infty
$$
with $\delta'>2\delta$
the norm $\V  A_N^{(p)}\V$ is bounded;
$$
{\hbox{p-}}\lim_{N,p\to\infty} \V A_N^{(p)}\V =2.
\eqno (2.3)
$$
}

\vskip 0.5cm
{\it Remarks.} 

1. In fact, we prove that under conditions of theorem 2.1
the estimate 
$$
\mu_a\otimes\mu_d \{\omega: \V A_N^{(p)}\V > 2(1+2\vep)\} \le 
\exp\{ - \varphi_N (\vep,\delta) \log N\}
\eqno (2.4)
$$
holds 
with  $\varphi_N(\vep\delta) = O((\log N )^\delta)$. This implies boundedness
of $\lim \V  A_N^{(p)}\V$ with probability 1.
 This fact together with (1.6) 
implies (2.3).

2. If there exists such $U$ that $V_{2k}\le U^{2k}$, then
(2.4) and  (2.3) hold for all $\delta'>0$.

\vskip 0.5cm

To show that 
condition $p\gg (\log N)^{1+\delta}$
 is necessary for convergence (2.3), we consider the simplest
case of Bernoulli random variables
$$
\hat a(x,y) = \cases{1,& with probability $1/2$;\cr
                -1,& with probability $1/2$.\cr}
$$

{\bf Theorem 2.2}

{\it Let $N,p$ increase infinitely in the way that
$$
{ p\over (\log N)^{1-\delta}}\to 0
\eqno (2.5)
$$
for any positive $\delta>1$.
Then for any given $R>0$
$$
\lim_{N,p\to \infty} \mu_a\otimes\mu_d \{\omega: \V \hat A_N^{(p)}\V > R \} = 1,
\eqno (2.6)
$$
where $\hat A_N^{(p)}$ is given by (1.4), (1.5) with 
$a(x,y)$ replaced by $\hat a(x,y)$.}
\vskip 0.5cm

Now we  describe the scheme of the proof of Theorem 2.1.
First of all, let us recall  that the semicircle law 
proved in \cite{W} states that
the moments
$$
M^{(N)}_j \equiv  \E \int \l^j \d \s (\l; \hat A_N) =\E {1\over N} \T \hat A_N^j 
\eqno (2.7)
$$
converges in the limit $N\to\infty$ to the numbers 
$$
m_j = \int \l^j\d\s_w(\l), \quad 
m_j = \cases{t_k,& for $j=2k$;\cr 
0,& for $j=2k+1$,\cr},
$$
where
$$
t_k= {(2k)!\over k! (k+1)!}.
\eqno (2.8)
$$

Regarding the average
$$
\E {1\over N} \T \hat A_N^j = 
{1\over N}\sum_{x=1}^N {1\over N^{j/2}} \sum_{y_i = 1}^N
\E a(x,y_1) a(y_1,y_2) \cdots a(y_{j-1},x)
\eqno (2.9)
$$
as a sum over paths 
$W(x,Y_{j-1})\equiv ( y_0 = x, y_1, y_2, \dots, y_{j-1},y_0)$ 
of $j$ steps, Wigner observed that
the leading contribution to (2.7) is given by the paths that
are related with the simple  walks in the upper half-plane
started and ended at zero. Basing on this fact, he derived
 recurrent relations for the limiting moments
$$
t_k = \sum_{l=0}^{k-1} t_{k-1-l} \, t_l, \quad t_0=1.
\eqno (2.10)
$$

In paper \cite{FK} the paths 
$W(x,Y_{j-1})$
were encoded with the help
of simple connected graphs with no cycles. 
More precisely, it was shown that the terms of
leading contribution to (2.9) can be related with the set $T_k$ of 
one rooted trees with
$k$ edges drawn in the upper half-plane. This representation is equivalent to
the simple walks description of \cite{W} and
the number $\v T_k\v$ of elements in $T_k$ is given by $t_k$.
Introducing some additional encoding for the trajectories $W$ 
that provide vanishing contribution,
the authors have studied asymptotic behaviour of the moments
$M_j^{(N)}$ for all $j\ll N^{1/6}, N\to\infty$
and proved boundedness of the limit of 
$\V \hat A_N\V$. These results were recently
improved in paper \cite{SS}, where 
the case of $j\ll N^{2/3}$ has been investigated.

Also let us note the paper \cite{BY}, where 
the necessary and sufficient condition for finiteness
of $\lim_{N\to\infty}\Vert A_N\Vert$ are found.
These conditions are rather weak and require 
existence of the fourth moment of the  random variables $a(x,y)$.

The methods developed in  \cite{BY} and \cite {SS} are  different from that used
in
\cite{FK}.  It is not clear whether they are applicable to the dilute random
matrices or not. 
In present paper we follow the way of \cite{W} and \cite{FK}
and develop it to be applied 
for the Wigner random matrices and their dilute versions as well.

There are three principal observations that we are based on:

(i)  all paths $W(x,Y_{2k-1})$ with non-zero contribution 
can be
separated uniquely into classes $\Pi(\tau)$ corresponding to the trees 
$\tau \in T_k$;

(ii) the paths that fall into the same class $\Pi(\tau)$ are described by
graphs $\g(\tau)$ that are obtained from $\tau$ by gluing its vertices;

(iii) the contribution of the path 
is given by the number of cycles  in $\g$. 

These facts allow one to estimate easily the number of graphs 
$\g$ and to sum the corresponding contributions.

In Section 3  we give precise account on the graph representation
of paths in the sum
$$
M^{(N,p)}_j \equiv \E \la N^{-1} \T  [A_N^{(p)}]^j\ra 
= 
{1\over N}{\sum_{x,y_i}}
\E A(x,Y_{2k-1}) \la D_N(x,Y_{2k-1})\ra
\eqno (2.11)
$$
where 
$$
 A(x,Y_{j-1})\equiv  a(x,y_1)a(y_1;y_2) \cdots a(y_{j-1},x)
$$
and
$$
D_N(x,Y_{2k-1})\equiv   \hat d_N(x,y_1) \cdots \hat
d_N(y_{j-1},x).
$$
It follows from this graphical description that
the terms of (2.11) that are of the order $O(p^{-s})$ 
are described by the graphs $\hat \g$ that 
have cycles only of the length 2. 
We show this in Section 4.

In Section 5 
we estimate the number of trees that can produce  
the cycles
of the length 2.
Basing on relation (2.10), we develop
a kind of the tree calculus and show that the number of trees
having vertices of large degree is exponentially
small with respect to the total number of trees.

In Section 6 we prove 
our main technical estimate
$$
M^{(N,p)}_{2k}\le  t_k \sum_{r,s=1}^k 
{(\alpha k)^{3r}\over  N^r} {(\beta k)^{s(1+\theta_k)}\over p^s} V_{4s}
\eqno (2.12)
$$
with $\theta_k\to 0$ as $k\to\infty$.
This inequality allows us to prove theorem 2.1. 
Theorem 2.2 is proved in Section 7.

We summarize our results and compare them
with certain known facts from probability theory and related fields
in Section 8. However, let us briefly discuss here one consequence of 
our results with respect to the Wigner random matrices 
that correspond to the case of $p=N$
in (1.5).
It follows from (2.12) with $p=N$ 
that one can assume that moments $\E \vert a(x,y)\vert^{2m}$ 
exists for $m\le m_0$ and  obtain boundedness of $M^{(N)}_{2k}$
for all $k\ll N^{\alpha}$, where $\alpha $ depends on $m_0$.
For example,  
using  a version of standard truncation procedure
(see e.g. a version presented in \cite{BY}),
one can easily show that
that $\alpha = 1/4$ requires $m_0=8$.
This indicates a bridge between results of \cite{BY} and \cite{SS}.
Relation (2.12) explains why the higher moments of 
the matrix entries $a(x,y)$
disappear from the limiting expressions for $M_{2k}$.

\section{Even partitions and trees}

It is clear that conditions (2.1a) imply that the average
$
\E A(x,Y_{j-1})
$
is non-zero if and only if the path 
$W(x,Y_{j-1})$ is even \cite{SS}, i.e. if each step $(y_i,y_{i+1})$ counted both 
in direct and inverse directions appears in $W$ even number of times.
This immediately lead to the conclusion that $M_{2k+1}^{(N,p)}=0$.

In the case when $j=2k$, an even path uniquely determines a partition $\pi$
of variables $\{y_i, i=0,\dots, 2k-1\}$ into groups. We call such a partition
the even one. 
The following statement links even partitions and rooted 
trees 
drawn in the upper half-plane.
Let us note that given such a tree, one can order its edges.
We adopt that the edges are enumerated  from below and from the left.
In other words, we say that the edge $e$ is less than the edge $e'$
when $e$ is situated on the path from the root to $e'$
or when there is a vertex $\nu$ 
such that the path from $\nu$ to $e$
is to the left with respect to the path
from $\nu$ to $e'$.

\vskip 0.5cm
{\bf Proposition 3.1.} 

{\it 
The set $\Pi'_{2k}$ of all even partitions of the variables
$(y_0,y_1,\dots y_{2k-1})$ can be uniquely separated into
classes $\Pi(\tau)$ labeled by the trees $\tau\in T_k$.
}

\vskip 0.2cm
{\it Proof.}  
We are going to show that 
given a partition $\pi\in\Pi'_{2k}$, 
one can  create a graph
$g(\pi)$ 
that is uniquely determined by a tree $\tau$.
We construct $g(\pi)$ with the help
of
natural procedure that
simplifies approach suggested in \cite{FK}
and resembles certain arguments of \cite{SS}.

The first step is to determine  the number of different
groups of variables in $\pi$.
This coincides with
the  number of vertices in $g(\pi)$. 
We refer to 
the  vertex that corresponds to the
group
containing $y_0$ as to {\it the root}.
Edges of $g(\pi)$ correspond to steps 
$s_i\equiv (y_i,y_{i+1})$, $i=0,\dots,2k-1$.

Next, one starts to go along $W=(y_0,\dots, y_{2k-1}, y_0)$ and
draws the edges of $g$ starting from the root and joining 
subsequently the vertices according to the order dictated by appearance
of steps $s_i$.  
One obtains closed connected graph with directed and ordered
edges.

Let us consider two vertices  joined
by one or more edges. Since all edges of $g$ are ordered,
the edges between two vertices are also ordered 
and can be numerated between themselves.
We call   the edge {\it marked} if it has an odd number 
in this inner enumeration.
It is not hard to see that there are 
exactly $k$ marked edges in $g$.

Finally one has to go once more along $g$  
and enumerate  the marked edges
 according to the order of their appearance in $g$. 
Regarding the marked edges only, one obtains the connected graph 
$\g$ with $k$ numerated edges
and one root. Let us call it {\it the first part} of the walk $g$.
The unmarked steps of $g$ make also connected graph with $k$ edges
that we call {\it the second part}  of $g$.

The graph $\g$ has one root and $k$ numerated edges. One can uniquely
restore the tree $\tau$ excepting the cases when
the edges $i$ and $i+1$ in $\g$ have both vertices common.
In this case we accept that $\tau$ is such that these vertices are consecutive.

 Proposition is proved.$\Box$

 Let us note in conclusion that 
 graph $\g(\tau)$ is obtained from $\tau$ by gluing  $q\ge 0$ vertices
 among them. 
One can regard this as  if $\g $ obtained by a procedure of $q$ steps where on
each step just one pair of vertices is glued.

Also it should be stressed that there exist partitions 
$\pi \neq \pi'$ that lead to the same graph $\g(\tau)$.
The difference between $\pi$ and $\pi'$
is that they
may have different
second parts of $g$ and $g'$.
This happens only when in $\g$
there are cycles of the length more than $2$
with edges "correctly" enumerated, i.e. one can go along the 
the cycle and see edges with increasing numbers 
prescribed for the edges of $\tau$.
We call these cycles as  {\it the correct loops}.
In all other cases the second part of $g$ is uniquely reconstructed from
$\g$.
We give precise account on this property  in Section 4.

\section{Number and contributions of partitions}

It is obvious that  
the value of $ \E A(x,Y_{2k-1})$ as well as 
that of $\la D_N(x,Y_{2k-1})\ra$
does not change if one preserves the partition and  moves variables $y_i$. 
Therefore one can write relation
$$
M^{(N,p)}_{2k}= \sum_{\pi\in \Pi'_k} B(\pi), \quad
B(\pi) \equiv {1\over N}{\sum_{x,y_i}}^{(\pi)}
\E A(x,Y_{2k-1}) \la D_N(x,Y_{2k-1})\ra
\eqno (4.1)
$$
and $\Sigma_{x,y_i}^{(\pi)}$ denotes the sum taken over variables 
$x,y_1, \dots,y_{2k-1}$ 
in the way such that the partition $\pi$ is preserved.

Tree representation helps to compute the contribution provided
by  the sum over given even partition $\Sigma^{(\pi)}$.
This contribution depends on the form of the graph $\g(\tau)$ 
that corresponds to $\pi$.

Let us introduce several terms. Assume that in a tree $\tau$
there is a vertex  $\nu$ of degree $m\ge 2$. 
The set of edges adjacent to $\nu$ we call {\it  the cluster}
of {\it the power} $m$. 
If two vertices of the edges belonging to the same cluster are glued,
then we say that there is {\it the cluster gluing} in the tree $\tau$.
The other gluings are called {\it the ordinary} ones.

\vskip 0.5cm 

{\bf Proposition 4.1}.

{\it 
Assume that 
$\g(\tau)$ is obtained 
from $\tau$ with the help of $q= r+s$ gluings, where $s$ is the
number of cluster gluings. 
Then 
$$
B(\pi) = B(\g) \le {1\over N^r\, p^s}\prod_{i=1}^j V_{4\xi_i}.
\eqno (4.2)
$$
}

{\it Remark. }
We say that the graph $\g$ with $r+s$ 
gluings and corresponding partition $\pi$
are {\it of the type $(r,s)$}.

\vskip 0.5cm
The proof of Proposition 4.1 can be easily  derived from
two observations.

The first one is that the number of vertices in $\g$ is equal to the number of
groups of variables in $\pi$. Thus
there are $(N-1)(N-2)\cdots(N-(k-q+1))$ terms in the sum $\Sigma^{(\pi)}$
of (4.1).

The second observation follows from the definitions  (1.4) and (1.5).
It is easy to see that if in $\g$ there exist 
two vertices $\mu$ and $\nu$ such that there are $\xi> 1$
edges $(\mu,\nu)$, then such a 
{\it multiple junction} provides the factor
$
\E a^{2\xi} \la \hat d_N^{2\xi}\ra = 
{V_{2\xi} N^{-1}p^{1- \xi}}.
$
If $\xi =1$, then the factor due to this simple
junction is obviously $1/N$.

It is clear that if there are $s$ gluings in one cluster,
then one obtains $l$ multiple junctions with multiplicities
$\xi_1, \dots ,\xi_l$ such that $\xi_1+\dots \xi_l = s+l$.
Taking into account elementary inequalities
$$
V_{2\xi_1}\cdots V_{2\xi_l}\le V_{2\xi_1+\dots 2\xi_l}=
V_{2(s+l)}\le V_{4s},
$$
one arrives at (4.2). $\Box$.

\vskip 0.5cm

Proposition 4.1 shows that all the partitions  of the type $(r,s)$ provide the
same contribution. 
Now it remains to estimate their number.
This number can be estimated by the number ${\cal N }(r,s)$ of possibilities
to make $r+s$ gluings in the tree $\tau$. We also should 
multiply ${\cal N }(r,s)$ by
the number of different partitions that correspond to the same
graph $\g$ of the type $(r,s)$.

Let us note here that the cluster gluings are independent from 
the other gluings in the sense that 
the number of possibilities to make the cluster gluing
in $\g $ does not depend on the number 
of other gluings already performed.
This observation
together with (4.2) shows that one can treat these two types of gluings
separately.

It is not hard to see that given a graph $\hat \g(\tau)$ of type $(o,s)$,
corresponding partition $\hat \pi$ can be uniquely restored.
This is because the way along $\hat \g$ preserving the existent order
can be performed
uniquely. In other words, the first and the second parts of the graph
$\hat g$ are uniquely determined by $\hat \g$.
Thus, the number of partitions of the type $(0,s)$
coincides with the number $D_{s}(\tau)$ of possibilities to make 
$s$ cluster gluings in $\tau$.

The picture differs for the ordinary gluings leading to the cycles
of the length greater than $2$. Assume that there is a vertex $\nu$
with $m>3$ edges $e_i$ and
there is a cycle of the length $l\ge 3$ starting and ending
at $\nu$. If this is the correct loop, then 
one can reconstruct the second part of the graph $g$ in several ways
because one can pass the cycle 
for the second time between passing the edges $e_i$.
Also each loop can be passed in two opposite directions and this
makes double the total number of possible ways.

Now let us estimate the number
of partitions of the type $(r,s)$
performed in a tree $\tau$ that  has 
$l$ clusters of powers $m_1,\dots ,m_l$,
respectively. Clearly, $k \le m_1+\dots m_l\le 2k-2$.
Let us assume that there are
$s$ cluster gluings such that in cluster $i$ there are $s_i$
gluings.
Then, summarizing the arguments presented above,
one can write that
$$
M^{(N,p)}_{2k}\le \sum_{\tau\in T_k} \sum_{q=0}^{k}\sum_{r+s=q} 
{1\over N^r}\left[ D_\tau (r) \sum_{\{R_l\}}^r \CN_\tau (\{R_l\})\right]
{1\over p^s}
\left[ \sum_{\{S_l\}}^s D_\tau(\{S_l\}) \prod_{i=1}^l V_{4s_i}\right].
\eqno (4.3)
$$
In this formula 
$D_\tau (r)$ is the number of possibilities
to make $r$ ordinary gluings in the tree $\tau$,
$\CN_\tau (\{R_l\})$ is the number of different partitions
that can be obtained from the graph $\g(\tau)$
due to presence of correct loops,
and $D_\tau(\{S_l\})$ is the number of possibilities
to make $s$ cluster gluings in $\tau$.
Here we have denoted 
$R_l = (r_1,r_2,\dots,r_l)$, 
and $S_l = (s_1,s_2,\dots,s_l)$,
such that $r_1+r_2+\dots+r_l=r$. Obviously, 
$s_1+s_2+\dots +s_l= s$.
Summations go over all possible combinations of $r_i\ge 0$
and $s_i\ge 0$ satisfying conditions presented above.

Elementary calculation shows that
$$
D_\tau (r) \le {1\over r!} {{k}\choose{2}}{k-1\choose 2} \cdots {k-r+1\choose 2}
\le {k^{2r}\over 2^r r!}.
\eqno (4.4)
$$
Next, assuming that all cycles obtained are correct, we can write that
$$
\CN_\tau(r)\equiv \sum_{\{R_l\}}^r \CN_\tau (\{R_l\} ) 
\le \sum_{\{R_l\}}^r \prod_{i=1}^l
{(m_i-2)!\over (m_i-1-r_i)!}.
\eqno (4.5)
$$
One can obtain the latter expression 
regarding the observation that if in a cluster $i$ there are $r_i$
correct loops, then the number of different ways to pass them
is estimated by the number of possibilities to distribute
$r_i$ different balls into $m_i-r_i$ boxes.
We easily derive from (4.5) that
$$
\CN_\tau (r) \le \sum_{\{R_l\}}^r  \prod_{i=1}^l 2^{r_i}{m_i^{r_i}}
\le
2^r\sum_{\{R_l\}}^r r! \prod_{i=1}^l {m_i^{r_i}\over r_i!}\le
(m_1+\dots+m_l)^r=(4k)^r.
\eqno (4.6)
$$

Let us turn to the number of cluster gluings.
It is estimated by the product
$$
\prod_{j=1}^s {1\over s_j!}{m_j\choose 2}{m_j-1\choose 2}\cdots{m_j-s_j+1\choose 2}\le
{1\over s!} {s!\over s_1!\cdots s_j!} \prod_{j=1}^s {m_j^{2s_j}}.
$$
Taking into account inequality
$
V_{4s_1}\cdots V_{4s_l}\le V_{4s},
$
we can write that
$$
\sum_{\{S_l\}}^s D_\tau(\{S_l\}) \prod_{i=1}^l V_{4s_i}\le
{ (m_1^2+\dots+m_l^2)^s\over s!} V_{4s}.
\eqno (4.7)
$$

Trivial inequality for the variable
$\Sigma_k^{(2)}(\tau) \equiv m_1^2 + \dots m_l^2\le 4k^2$ lead to
the estimate
$k^{2s}/p^s$.
This estimate is rather rough and  is not sufficient for the proof of theorem 2.1.
In the next section we show that 
$\sum_\tau \Sigma_k^{(2)}(\tau)$  behaves like
$t_k k^{1+\delta}$ that  reflects the almost linear
character of the average tree. This fact  together
with inequalities (4.5)-(4.7) implies (2.12).

\section {Enumeration of trees}

In this section we  consider the set $T_k$ the one-root
trees $\tau$ having $k$ edges drawn in the upper half-plane. 
Let us recall that given a tree $\tau\in T_k$, one can
order its edges. 
Let us also give several  definitions.

We refer to the edges adjacent to the root as to the {\it root edges. }
We define the {\it cluster} $\kappa $ as the set of $m$ edges that have one
common vertex $\nu$. We call $m$ and $\nu$ the {\it power} 
and the center, respectively, of the cluster $\kappa $.

The main result of this section is given by the following statement.

\vskip 0.5cm

{\bf Theorem 5.1.}

{\it Let us consider the set $G_k^{(m)}\subset T_k$ of trees
that have at least one cluster with the power $m\ge 2$. Then the number
of such trees $|G_k^{(m)}|\equiv g_k(m)$ satisfies relation
$$
g_k(m)\le  k t_{k}\left( \frac 34\right) ^{m-2}.
\eqno (5.1) 
$$
}

To prove theorem 5.1, we need the following auxiliary statement.
\vskip 0.5cm

{\bf Lemma 5.1.}

{\it Let us consider the set $T_k^{(m)}$ of trees that have $m\ge 2$ root
edges. Then} 
$$
t_k^{(m)}\equiv |T_k^{(m)}|\le t_{k-1}\left( \frac 34\right) ^{m-2}.
\eqno (5.2) 
$$

{\it Proof. } 
It is easy to see that 
$$
t_k^{(m)}=\sum_{\alpha _i}^{k-m}\, t_{\alpha _1} t_{\alpha _2}\cdots
t_{\alpha _m},
\eqno (5.3) 
$$
where the sum is taken over all $\alpha _i\ge 0$ such that $\alpha _1+\dots
+\alpha _m=k-m$. 

Let us derive first  the simple estimate 
$$
t_k^{(m)}\le t_{k-1}.
\eqno (5.4) 
$$
Regarding (2.10), we can rewrite (5.3) in the form 
$$
t_k^{(m)} =\sum_{q=0}^{k-m}
\sum_{\alpha _i}^{k-m-q}\,t_{\alpha _1}t
_{\alpha _2}\cdots t_{\alpha _{m-2}}\sum_{\alpha _i}^q\,t_{\alpha_{m-1}} 
t_{\alpha _m}= 
$$
$$
\sum_{q=0}^{k-m}
\sum_{\alpha _i}^{k-m-q}\,t_{\alpha _1}t
_{\alpha _2}\cdots t_{\alpha _{m-2}} t_{q+1}.
\eqno (5.5)
$$

The latter sum can be regarded as the sum over $m-1$ variables $\a_i\ge 0$,
where the term with $\a_{m-1}=0$ is absent. Then
$$
t_k^{(m)}=\sum_{\alpha _i}^{k-m+1}\, t_{\alpha _1} t_{\alpha _2}\cdots
t_{\alpha _{m-1}}-
\sum_{\alpha _i}^{k-m+1}\, t_{\alpha _1} t_{\alpha _2}\cdots
t_{\alpha _{m-2}} \equiv 
t_k^{(m-1)} - t_k^{(m-2)}.
$$
Thus, $t_k^{(m)}\le t_k^{(2)}$. Relation (2.10) 
implies equality $t_k^{(2)} =t_{k-1}$ that gives (4.4).

Let us rewrite (5.3) in the form
$$
t_k^{(m)}= \sum_{q=0}^{k-m} t_{k-q-1}^{(m-1)} t_q.
\eqno (5.6)
$$
If $m-1\ge 2$, then we can apply (4.4) to 
$t_{k-q-1}^{(m-1)}$ and obtain that
$$
t_{k}^{(m)}\le
\sum_{q=0}^{k-m} t_{k-q-2} t_q=
$$
$$
t_{k-1} - (t_0t_{k-2} + \dots + t_{m-1}t_{k-m-1})\le t_{k-1}-t_{k-2}.
$$
Expression (2.8) for $t_k$ implies that $t_{k-2} > t_{k-1}/4$.
Therefore $t_{k}^{(m)}\le 3t_{k-1}/4$. 

If $m-2$ is greater than $2$, then we can substitute the last inequality
into (5.6) and obtain that 
$t_{k}^{(m)}\le t_{k-1}(3/4)^2$. Now it is clear that (5.2) is true.

Using similar computations, one can easily prove the following statement.
\vskip 0.5cm
{\bf Lemma 5.2}.

{\it Let us denote by $n_r(s)$ the number of trees that can be constructed 
on $r\ge 1$ roots with the help of $s$ edges.
Then for $2\le l\le s$}
$$
n_{r+l}(s-l) \le n_{r+1}(s-1)  \left( {3\over 4}\right)^{l-2}\le
 n_{r}(s) \left( {3\over 4}\right)^{l-1}.
\eqno (5.7)
$$

{\it Proof theorem 5.1.} 

We start with the  observation
that one can construct 
a tree $\tau\in T_k$
from the set $E_k$ of $k$ edges $e_i$ that are already enumerated 
on the way that 
this enumeration will
agree with the order among tree edges.

Now assume that  before this  procedure the edges 
$e_h, e_{h+1},\dots, e_{h+m}$  are joined to the same vertex 
and used in the construction as one cluster. According to Lemma 4.2,
this diminishes exponentially with respect to $t_k$ 
the number of trees obtained on this way.
Taking into account that $h$ can vary from $1$ to $k-m$
we obtain the estimate (5.1).$\Box$
\vskip 0.5cm

In conclusion let us note that  (5.1) implies that the number
of trees that have the power of the maximal cluster 
large enough is exponentially small with respect to $t_k$.
Therefore the number of such clusters is exponentially small
with the total number of clusters in trees $\tau\in T_k$.
This means that the average degree of a tree vertex remains
finite even in the limit $k\to\infty$.
We plan to study this problem mores systematically
in a separate publication.

\section {Proof of theorem 2.1.}

As it follows from (5.3) and (5.7), we need to estimate the sum
$$
Q_k(s) \equiv \sum_{s=1}^{k-1} {V_{4s}\over s! p^s} 
\sum_{\tau\in T_k} 
(m_1^2+\dots +m_l^2)^s,
$$
where $m_i$ are the powers of clusters of $\tau$. 
Let us denote by $\hat m(\tau)$
the maximal cluster power of the tree $\tau$; then obviously
$$
\sum_{\tau\in T_k} 
(m_1^2+\dots +m_l^2)^s \le \sum_{\tau\in T_k} 
\hat m(\tau)^s (m_1+\dots +m_l)^s = (2k)^s 
\sum_{\tau\in T_k} 
\hat m(\tau)^s.
$$

Let us consider the set of trees $T_k'$ such that
$\hat m_\tau \le \chi \log k$, where $\chi = (\log 4 - \log 3)^{-1}$.
Then
$$
(2k)^s 
\sum_{\tau\in T_k'} 
\hat m_\tau^s \le t_k (2k \chi \log k)^s \le t_k (2\chi k)^{s(1+\theta)}
$$
with some $\theta=\delta_k$ that can be taken vanishing as $k\to\infty$.

We represent the set $T_k\setminus T_k'$ as the sum of sets
$G_k( [\chi(j+\log k )] )$ of trees that have the maximal cluster
power equal to $[\chi(j+\log k)]$, $j=1,2,\dots$.
Here we denoted by $[x]$ the maximal natural number less than 
or equal to $x$. 
Using estimate (5.1), we can write that
$$
\sum_{\tau \in T_k'} \hat m_\tau^s\le
\sum_{j=1}^{k-[\chi \log k]} 
\sum_{\tau \in G_k(j+[\chi\log k] ) } \hat m_\tau^s\le
$$
$$
e^2 \sum_{j=1}^k \chi^s (\log k +j)^s \exp\{ - j - \log k\} k t_k \le
e^2 (s+1)! (\chi \log k)^s.
$$
Thus, we obtain inequality (cf. (2.12)) 
$$
Q_k(s)\le p^{-s} (\beta k)^{s(1+\theta')} V_{4s},
\eqno (6.1)
$$
where $\beta$ is a constant and 
$\theta'$ can be chosen vanishing when $k\to\infty$.

Now it is easy  to derive the estimate (2.4) from (2.12).
Taking into account (6.1) and the estimate $ V_{4k}\le (2k)^{2\delta k}$,
we obtain inequality
$$
M_{2k}^{(N,p)}\le t_k 
\left( (1 - \alpha k^3N^{-1}) (\beta k^{1+2\delta
+\theta'} p^{-1})
\right)^{-1}.
$$
Then we can deduce  that inequality
$$
M_{2k}^{(N,p)}\le \vep t_k\le (1+\vep)^{2k}\,4^k
\eqno (6.2)
$$
holds for all $k$ sufficiently large and such that $k\ll N^{1/3}$ and 
$k^{1+\delta+ \theta'}\ll p$.

Regarding definition (1.2), we can write that 
$$
M_{2k}^{(N,p)}\le \E\la \int_{\vert\l\vert \ge 2(1+2\vep)} \l^{2k} \d
\s(\l;A_N^{(p)})\ra \ge
4^k (1+2\vep)^{2k} N 
\E\la \#\{\vert \l_j^{(N,p)} \vert \ge 2(1+\vep)    \} \ra\ge
$$
$$
4^k (1+2\vep)^{2k} N \mu_a\otimes\mu_d \{\omega: \V A_N^{(p)}\V > 2(1+2\vep)\}.
$$
Combining these inequalities with (6.2), we derive that
$$
\mu_a\otimes\mu_d \{ \omega: \V A_N^{(p)}\V > 2(1+2\vep)\} \le 
N \left( {1+\vep\over 1+2\vep}\right)^{2k}.
$$
Choosing $k$ such that $k/\log N\to\infty$, we obtain 
that (2.4) holds. Theorem 2.1 is proved.

\section {Proof of theorem 2.2}

Let us consider the unit vectors $\vec e^{(j)}$ with the components 
$$
\vec e^{(j)}(x)=\cases{1,& if $x=j$,\cr
0,& if $x\neq j$.\cr} 
$$
We can write that%
$$
\left\| A_{N,p}\right\| ^2\ge \max {}_{j=1,\dots ,N}\;\left\| A_{N,p}\vec
e^{(j)}\right\| ^2=\max {}_{j=1,\dots ,N}\left( A_{N,p}^2\right) (j,j), 
\eqno (7.1) 
$$
where%
$$
\left( A_{N,p}^2\right) (j,j)=\sum_{y=1}^N\left| A_{N,p}(j,y)\right| ^2. 
$$
Let us introduce random variables
$$
h_j^{(N)}=\sum_{y=j}^N\left| A_{N,p}(j,y)\right| ^2,\;j=1,\dots ,N. 
$$
It is clear that the family $\left\{ h_j^{(N)}\right\} _{j=1}^N$ is the set
of jointly independent random variables and that
$$
\max_{j=1,\dots ,N}\left( A_{N,p}^2\right) (j,j)\ge \max_{j=1,\dots
,N}\;h_j^{(N)}.
\eqno (7.2) 
$$
Let us note that under conditions of theorem 2.2 
$$
h_j^{(N)}=\frac 1p\sum_{y=j}^N d(j,y)\equiv \frac 1p\eta _j^{(N)}.
\eqno (7.3)
$$

Let us consider the probability distribution of the random variable 
$$
H^{(N)}=\max_{j\le N/2}\;h_j^{(N)}. 
$$
It is clear that 
$$
P_N(R)\equiv \Pr \left\{ H_N<R \right\} =\Pr \left\{ \max_{j\le N/2}\;\eta
_j^{(N)}<pR\right\}  
$$
$$
=\prod_{j=1}^{N/2}\left( 1-\Pr \left\{ \eta _j^{(N)}\ge pR\right\} \right) . 
$$
It is not hard to prove that the random variables $\eta _j^{(N)}$ considered
for $p\sim \log N$ converge as $N\to\infty$ 
to the  Poissonian random variables
$\zeta_j$  that are identically distributed with the parameter $p$. Then we can
write that
$$
1-\Pr \left\{ \eta _j^{(N)}\ge pR\right\} \le 1-\Pr 
\left\{ \zeta_j=pR\right\}
=1-\frac{p^{pR}}{(pR)!}e^{-p}. 
$$
Let us show that if $p=\left( \log \;N\right) ^{1-\delta },$ then 
$
p^{pR}\left[ e^p(pR)!\right] ^{-1}$ decays more slowly than $2/N$ and
$$
\left[ 1-\frac{p^{pR}}{(pR)!}e^{-p}\right] ^{\frac N2}\rightarrow 0\;{as
}\;N\rightarrow \infty .
\eqno\,(7.4) 
$$
Using the Stirling formula, we can write that
$$
e^{-p}\frac{p^{pR}}{(pR)!}\sim \frac{e^{p(R-1)}}{R^{pR}}=\exp \left\{
-pR\left[ \log R-\frac r{r-1}\right] \right\} . 
$$
It is easy to see that (7.4) holds that together with (7.1) and (7.2) proves
relation (2.6).

\section {Summary and discussion}

We consider the ensemble of random matrices
$A_{N,p}$ with independent entries 
such that $A_{N,p}$ has, in the average
$p$ non-zero entries per row.
We study asymptotic behaviour of the spectral norm $\V A_{N,p}\V$
in the limit when $N,p\to\infty$ and $p=o(N)$. We consider the averaged moments
$$
M_{2k}^{(N,p)}= {1\over N}  \sum_{x=1}^N \sum_{\{y_i\}} \E \la
A_{N,p}(x,y_1)\cdots A_{N,p}(y_{2k-1},x)\ra
\eqno (8.1)
$$ 
and give a further development to the approach proposed in papers
\cite{FK,W}.
The method is based on the relation between   the even partitions $\pi$
of variables $(x,y_1, \dots,y_{2k-1})$, i.e. those that
make the average in (8.1) non-zero, from one hand and the rooted trees $\tau\in T_k$
constructed in the upper half-plane with the help of $k$ edges, from the other hand.

We give the full description for the even partitions $\pi$ in terms of the 
graphs $\g$ obtained from a tree $\tau$ by gluing its vertices.  
The magnitude of the terms described by  partition $\pi$ 
is determined by the number of cycles
in corresponding graph $\g$. This allows one to estimate
easily the number of different partitions and their 
contributions to $M_{2k}^{(N,p)}$.

We show that the terms of the order $p^{-l}$
$M_{2k}^{(N,p)}$ in the limit $p=o(N),\,N\to\infty$ 
comes from those partitions $\pi'$
that are encoded by the graphs $\g'$ that have cycles of the length 2 exactly.
This means that the average degree of the vertex plays
the key role in asymptotic properties of large sparse random matrices.

We obtain the estimate showing that the average in certain sense tree
is of the linear structure and the relative number of possible 2-cycles increases
almost linearly with respect to $k\to\infty$. 
This allows us to derive our main result  that $\Vert A_{N,p}\Vert$
remains bounded when $p\sim (\log N)^{1+\delta}, \,\delta>0$. 
We show that the 
value $p_c=\log N$ is  the critical one because
$\Vert A_{N,p}\Vert$ is unbounded when $p\sim (\log N)^{1-\delta}$.

Our results can be compared with the classical Erd\"os-R\'enyi limit
theorem \cite{ER}. This statement concerns independent random variables
$\xi_1,\dots, \xi_N$ and 
the partial sums $X_{i,p}= \xi_i+\xi_{i+1}+\dots+\xi_{i+p}$.
It is proved that  the value $p_c=\log N$ is the critical one
for the $X(N)=\max_i X_{i,p} p^{-1}$ to be either bounded or not.

One can regard our theorems 2.1 and 2.2
as a limit theorems for stochastic versions of Erd\"os-R\'enyi partial sums.
This claim is supported by the fact that random variables
$h_j^{(N)}$ (5.3) are determined as the sum of approximately
$p$ random variables.
Let us also note that theorem 2.1 concerns a random variable
$\Vert A_{N,p}\Vert$
that is greater  than $\hat X(N)\equiv \max_j \Vert A_{N,p} \vec e^{(j)}\Vert$.
This could explain the need of conditions (2.2) that are more restrictive
than the standard conditions of Erd\"os-R\'enyi limit theorem.

One can also trace out more subtle link between our results and
the well-known statement concerning the connectedness  of
random graphs (see e.g.\cite{B}). If one has $N$ vertices and draws
$q$ edges joining randomly chosen vertices, then the graph obtained
will be connected in the limit $N\to\infty$
provided $q\gg N\log N$.
In our terms the randomly chosen pairs $(i,j)$ 
to be joined correspond to non-zero elements in the dilute
matrix $A_{N,p}$. In this context, the  transition from disconnected 
graph to the connected 
one can be regarded as analog of the transition of   
$A_{N,p}$  from the class of the tridiagonal
matrices (like discrete Schr\"odinger operator)
to the class of Wigner random matrices.

It should be noted that 
in mathematical physics literature another critical values 
are  found for sparse random matrices.
In particular, in paper \cite{MF} it is claimed that 
certain "density-density" correlation function changes its behaviour
at finite values of $p$. Numerical studies \cite{E} show that 
certain spectral characteristics of strongly dilute random matrices 
can depend on finite values of $p$.
Therefore our results imply that there can be several
different critical values in the sparse random matrix model.

To complete the discussion, let us note that the dilute random matrices
$A_{N,p}$ with $p$ replaced by $N$ take the form of the Wigner random
matrices $\hat A_N$. Therefore technique developed in present paper
can be also useful in this case. This 
topic is already  addressed  in Section 2.

\vskip 0.5cm
{\bf Acknowledgments.} The author is grateful to Dr. V.Gayrard and
Dr. A.Soshnikov for useful discussions and to Prof. F.Comets
for his comments about the Erd\"os-R\'enyi limit theorem.
The author also thanks Laboratory of Mathematical Physics and Geometry,
University Paris-7, where the paper is completed, for kind hospitality
and Prof. A.Boutet de Monvel for her interest to this work.

\end{large}
\end{document}